\documentclass[3p,twocolumn,number]{elsarticle}
\usepackage{paralist}
\usepackage{amsmath,amssymb,color,tikz}
\usepackage[numbers]{natbib}
\newtheorem{theorem}{Theorem}

\newtheorem{conjecture}{Conjecture}

\newcommand{\beq}{\begin{equation}}
\newcommand{\eeq}{\end{equation}}

\newcommand{\E}{\mathbb{E}}
\renewcommand{\P}{\mathbb{P}}
\newcommand{\V}{\mathbb{V}\mbox{ar}}
\renewcommand{\O}{\mathcal{O}}

\newcommand{\LST}[1]{\widetilde{#1}}
\newcommand{\SCV}[1]{\textit{cv}^2_{#1}}

\begin{document}

\title{Critically loaded k-limited polling systems}

\author[tue]{M.A.A. Boon\corref{cor1}}
\ead{m.a.a.boon@tue.nl}
\cortext[cor1]{Corresponding author}
\address[tue]{Eurandom and Department of Mathematics and Computer Science, Eindhoven University of Technology, P.O. Box 513, 5600MB Eindhoven, The Netherlands.}

\author[uva]{E.M.M. Winands }
\ead{e.m.m.winands@uva.nl}
\address[uva]{University of Amsterdam, Korteweg-de Vries Institute for Mathematics, Science Park 904, 1098 XH  Amsterdam, The Netherlands.}

\begin{abstract}
We consider a two-queue polling model with  switch-over times and $k$-limited service (serve at most $k_i$ customers during one visit period to queue $i$) in each queue. The major benefit of the $k$-limited service discipline is that it - besides bounding the cycle time - effectuates prioritization by assigning different service limits to different queues. System performance is studied in the heavy-traffic regime, in which one of the queues becomes critically loaded with the other queue remaining stable. By using a singular-perturbation technique, we rigorously prove heavy-traffic limits for the joint queue-length distribution. Moreover, it is observed that an interchange exists among the first two moments in service and switch-over times such that the HT limits  remain unchanged. Not only do the rigorously proven results readily carry over to  $N$($\geq2$) queue polling systems, but one can also easily relax the distributional assumptions. The results and insights of this note prove their worth in the performance analysis of Wireless Personal Area Networks (WPAN) and mobile networks, where different users compete for access to
the shared scarce resources.
\end{abstract}


\maketitle

\section{Introduction}\label{sect:introduction}

This note studies a two-queue $k$-limited  polling model  in heavy-traffic (HT), in which one of the queues becomes critically loaded with the other queue remaining stable. Under the $k$-limited strategy the server continues working  until either a predefined number of $k_i$ customers is served at queue~$i$ or until the queue becomes empty, whichever occurs first. The $k$-limited policy is easy to implement, augmenting the exhaustive strategy (where $k_i=\infty$) by a single ``knob'' for each queue that can tune performance to be efficient or fair.

Significant gains in system performance have been observed in many application areas, e.g. communication \cite{borst,charzinski}) and production systems (see, e.g. \cite{winands}), by implementing the $k$-limited strategy instead of the traditional exhaustive and gated strategies (see, also, \cite{boonapplications2011} for a survey of applications of polling systems). Our specific interest for this model is fueled by performance analysis of (Bluetooth) Wireless Personal Area Networks (WPAN) and mobile networks, where different users compete for access to the shared scarce resources. In case of a Round-Robin type of scheduling, polling models with limited service policies naturally emerge (see \cite{boonapplications2011}). Proper operation of these systems and the scheduling therein is particularly critical when the systems are critically loaded.

A mathematical challenge is that  the $k$-limited strategy does not satisfy the branching property for polling systems, which significantly increases the analytical complexity \cite{resing}. A testament for this statement is that - notwithstanding the wealth of literature on polling systems - the problem of determining the steady-state distribution for our model is one that remains open. The statement is reinforced as well by the fact that heavy-traffic limits have only been rigorously proven for  branching-type service disciplines with very few exceptions. Non-branching policies differ from their branching counterparts both in the techniques useable and the resulting HT behavior. That is, the techniques used for branching-type policies typically rely heavily on explicit steady-state results for stable systems; however, corresponding results for non-branching
policies are not available. Moreover, non-branching systems can possess both stable and instable queues which is very different from branching policies where all queues become instable simultaneously leading to incomparable HT behavior.

We rigorously prove HT asymptotics for two-queue $k$-limited polling models with switch-over times via the singular-perturbation technique. That is, by increasing the arrival intensity at one of the queues until it becomes critically loaded, we derive HT limits for the joint queue-length distribution. The singular-perturbation technique can be typified by its intrinsic simplicity and its intuitively appealing derivation, although it requires some distributional assumptions (see \cite{knessltiersurvey} for a survey of applications of the
perturbation technique to queueing models). This technique enables us to present the first study rigorously proving HT asymptotics for polling models with a non-branching service strategy. It is proven that the scaled queue-length distribution of the critically loaded queue is exponentially distributed, whereas the number of customers in the stable queue has the same distribution as the number of customers in a vacation system. Moreover, both queue-length processes are independent in HT. The analysis can be directly extended to an $N$-queue system ($N \geq 2$) with one queue becoming critically loaded, as well as one can relax the distributional assumptions.

As a by-product, we observe a close  similarity between $k$-limited polling systems with and without switch-over times. This novel result, which is reminiscent of the variance absorption result in branching-type polling systems \cite{cooper}, shows that in HT the first two moments of the switch-over times can be effectively absorbed into the corresponding moments of the service times without impacting the HT asymptotics of the instable queue. In this regard, new fundamental insight is gained which should prove useful in the analysis of non-branching polling systems.

 The approach in the present note has its origin in \cite{boonperturbation}, where systems with \textit{zero} switch-over times  are studied. At face value the extension to \textit{nonzero} switch-over times may seem a small one, however this extension impels us to, considerably, modify and extend the analysis in \cite{boonperturbation}. This reveals itself clearly in the determination of the parameter of the exponential distribution for the scaled queue length of the critically loaded queue. In case of zero switch-over times this parameter follows directly from the HT limit of the total workload in a standard $M/G/1$ model and the averaging principle.  When analyzing a system with nonzero switch-over times, service capacity is obviously lost due to the switching and, as a consequence, one cannot directly state that the total workload equals the  $M/G/1$ workload and the aforementioned interchange of switch-over and service time moments is required.  Lastly, we note that in many practical applications the switch-over times may be substantial and that
the presence of these switch-over times may be crucial for the operation of the system (see, e.g., \cite{winands}). In the interest of space, we do not give all proofs in the current note, but we focus on the main theorems, the theoretical insights resulting from these theorems and the observed similarity between polling systems with and without switch-over times.

\section{Model description}

We consider a polling model consisting of two queues, $Q_1$ and $Q_2$, that are alternately visited by a single server.  When $k_i$ customers have been served or $Q_i$ becomes empty, whichever occurs first, the server switches to the other queue. Then, an exponentially distributed switch-over time $S_i$ with parameter $\sigma_i$ is incurred.  If the other queue turns out to be empty, the server switches back and serves, again, at most $k_i$ customers. The situation where both queues are empty will turn out to be negligible in the HT limit, which is considered in this note. A possible scenario is that the server waits until the first arrival and switches to the corresponding queue (say, $Q_j$) to start another visit period of at most $k_j$ customers. Customers arrive at $Q_i$ according to a Poisson process with intensity $\lambda_i$. We assume that the service times of customers in $Q_i$ are independent and exponentially distributed with parameter $\mu_i$. We denote the load of the system by $\rho = \rho_1+\rho_2$, where $\rho_i=\lambda_i/\mu_i$. The utilization $u_i$ for queue $i$, $i=1,2$, is defined as follows
\begin{equation}
u_i = \rho + \lambda_i\frac{\E[S]}{k_i},\label{eqn:ui}
\end{equation}
where $S$ is the total switch-over time in a cycle, i.e. $\E[S]={1}/{\sigma_1}+{1}/{\sigma_2}$. The cycle time $C$ is defined as the time between two successive visit beginnings to $Q_i$, with mean $\E[S]/(1-\rho)$. A known result for this model (cf. \cite{fricker}) is that a necessary and sufficient stability condition for queue $i$ reads $u_i < 1$, $i=1,2$ \cite{fricker}. This condition can be rewritten  as follows,
\begin{equation}
\lambda_i \mathbb{E}[C] < k_i.
\end{equation}
In words, this means that for a stable system the average number of type-$i$ customers arriving in a cycle is smaller than the service limit $k_i$, i.e., the maximum number of type-$i$ customers served in a cycle. Note that this system might have one stable queue (with utilization less than 1) and one instable queue (with utilization greater than 1), in contrast to the commonly studied \emph{branching-type} polling systems where all queues become instable simultaneously (see, for example, Resing \cite{resing}).

\section{Analysis}

We study the limiting behavior of the model in HT. In particular, we take the limit such that the utilization of $Q_2$ reaches its critical limit 1, while $Q_1$ remains stable. In this note we take the limit by increasing the arrival rate of $Q_2$ until it reaches its critical value while keeping the other parameters fixed.
Note that increasing $\lambda_2$ will also increase $u_1$, the utilization at the other queue. In order to prevent $Q_1$ from becoming instable, we impose the following condition:
\[\lambda_1/k_1<\lambda_2/k_2.\]
The results still hold when other (combinations of) parameters are modified such that $u_2\uparrow1$. See Section \ref{discussions} for a short discussion on this topic. The case where both queues become instable simultaneously requires a different analysis leading to different results and is beyond the scope of this note (see Section \ref{discussions}).

\subsection{Perturbed balance equations}

The model is a continuous-time Markov chain, and we start by describing its states $(n_1,n_2,h)$, where $n_1$ and $n_2$ denote the queue lengths. The possible server states are numbered $h=1,2,\dots,k_1+k_2+2$ and should be interpreted as follows. States $1\leq h\leq k_1$ indicate that the server is serving the $h$-th customer during its current visit to $Q_1$. State $h=k_1+1$ indicates that the server is switching from $Q_1$ to $Q_2$. Similarly, values $k_1+2\leq h\leq k_1+1+k_2$ indicate that the server is serving the $h-k_1-1$-th customer at $Q_2$. The final state, $h=k_1+k_2+2$, indicates that the server is switching from $Q_2$ back to $Q_1$. We refrain from giving all transition rates as they are trivial to determine.

The next step, after having defined the states and transition rates, is to give all balance equations so we can apply a singular-pertubation to them. We select \emph{one} balance equation to illustrate the technique,
\begin{multline}
(\lambda_1+\lambda_2+\sigma_1)p(0, n_2, k_1+1) = \lambda_2p(0,n_2-1,k_1+1) \\
+\mu_1\sum_{h=1}^{k_1}p(1,n_2,h)+\sigma_2p(0, n_2, k_1+k_2+2),
\label{balanceeqn}
\end{multline}
for $n_2>1$. This equation relates the in-flow and out-flow of state $(0, n_2, k_1+1)$, which is the state where $Q_1$ is empty, $Q_2$ has $n_2$ customers, and the server is switching from $Q_1$ to $Q_2$. One can leave this state through an arrival at any of the queues, or because the switch-over time has ended. This state can be entered from state $(0,n_2-1,k_1+1)$, through an arrival at $Q_2$, or through a departure from $Q_1$ when in any of the states $(1,n_2,h)$ for $h=1,\dots,k_1$, or from state $(0, n_2, k_1+k_2+2)$. In the last case the switch-over time $S_2$ has finished, but since $Q_1$ is empty, the server starts switching back to $Q_2$ again.

We now introduce the following definitions
\begin{equation}
\E[B_2']=\frac{1}{\mu_2}+\frac{\E[S]}{k_2},\quad
\V[B_2']=\frac{1}{\mu_2^2}+\frac{\V[S]}{k_2},\label{meanvarB2}
\end{equation}
where $\V[S]=\frac{1}{\sigma_1^2}+\frac{1}{\sigma_2^2}$. The random variable $B_2'$ will be given an interpretation   in Section {\ref{discussions}.\\
\ \\
\textbf{Perturbation.} We increase the arrival rate of $Q_2$ as follows to its critical value,
\begin{equation}
\lambda_2=\frac{1-\rho_1}{\E[B_2']}-\delta\omega, \qquad \omega>0,0<\delta\ll 1.
\label{perturbation}
\end{equation}
Substituting \eqref{perturbation} in \eqref{eqn:ui}, taking $i=2$, and letting $\delta\downarrow0$ implies that $u_2\uparrow1$ and, thus, that $Q_2$ becomes instable. An appropriate value for the constant $\omega$ is chosen at the end of this section.

Let $\xi=\delta n_2$, and
\begin{equation}
p(n_1, \xi/\delta, h) = \delta \phi_{n_1, h}(\xi,\delta), \label{phi}
\end{equation}
for $0<\xi=\O(1), h=1,2,\dots,k_1+k_2+2$. Note that once $\omega$ is chosen, $\delta$ and the scaled variable $\xi$ will be uniquely defined.
The next step is to substitute \eqref{perturbation} and \eqref{phi} in the balance equations. For Equation \eqref{balanceeqn} we obtain,
\begin{multline*}
(\lambda_1+\sigma_1)\phi_{0, k_1+1}(\xi,\delta) - \mu_1\sum_{h=1}^{k_1}\phi_{1,h}(\xi, \delta)
- \sigma_2\phi_{0,k_1+k_2+2}(\xi,\delta)\\
= \left(\frac{1-{\lambda_1}/{\mu_1}}{\E[B_2']}-\delta\omega\right)\big(\phi_{0, k_1+1}(\xi-\delta,\delta)-\phi_{0, k_1+1}(\xi,\delta)\big).
\end{multline*}
Taking the Taylor series with respect to $\delta$ yields a set of equations.

\subsection{Main results}\label{mainresults}
In this subsection, we derive the limiting joint distribution of $N_1$ and $\delta N_2$, the queue length at $Q_1$ and the scaled queue length at $Q_2$, as $\delta\downarrow0$. The three results are obtained by equating the $\O(1)$, $\O(\delta)$, and $\O(\delta^2)$ terms in the set of equations from the previous subsection. These steps follow the line of reasoning introduced in \cite{boonperturbation} and, therefore, the details are omitted in interest of space.

We define $\LST{L}^{(0)}(z)$ as the steady-state queue-length PGF of a single-server multiple-vacation queue with parameters $\lambda_1$ and $\mu_1$, with $k_1$-limited service. The vacation distribution in this system equals the convolutions of $k_2$ service-time distributions in $Q_2$ and of the switch-over times $S_1$ and $S_2$. We also define $P_0(\xi)$ as the scaled queue-length distribution of $Q_2$, which is yet to be determined. Equating the $\O(1)$ terms in the perturbed balance equations gives the following result.
\begin{theorem} \label{theorem1}
\begin{equation}
\sum_{n_1=0}^\infty\sum_{h=1}^{k_1+k_2+2}\phi^{(0)}_{n_1,h}(\xi)z^{n_1} = \LST{L}^{(0)}(z)P_0(\xi).\label{solutionOrder0}
\end{equation}
\end{theorem}
For this vacation system, the PGF  of the queue-length distribution $\LST{L}^{(0)}(z)$ is known (cf. \cite{lee89}).

Subsequently, we determine the unknown  $P_0(\xi)$ in (\ref{solutionOrder0}) by solving the equations for the first-order and second-order terms in the perturbed balance equations. This results in
\begin{equation}
C_0 P_0''(\xi)+C_1 P_0''(\xi)+C_2 P_0'(\xi)=0,\label{differentialeqn}
\end{equation}
where the constants $C_0$, $C_1$, $C_2$ are defined as follows,
\begin{align*}
C_0 &= \frac{2\frac{\lambda_1}{\mu_1^2}+\frac{1}{\E[B_2']}\left(1-\frac{\lambda_1}{\mu_1}\right)\left[\E[B_2'^2]-2\E[B_2']^2\right]}{2\E[B_2']^2},\\
C_1 &= \frac{1}{\E[B_2']}\left(1-\frac{\lambda_1}{\mu_1}\right), \quad C_2 = \omega .
\end{align*}
%
%
%

We obtain the density $P_0(\xi)$ of the scaled number of customers in $Q_2$ by solving the differential equation \eqref{differentialeqn}, where we use that $\sum_{n_1=0}^\infty\sum_{n_2=0}^\infty\sum_{h=1}^{k_1+k_2}p(n_1,n_2,h)=1$ and that $\int_0^\infty P_0(\xi)\textrm{d}\xi=1$.

\begin{theorem} \label{theorem2}
\begin{equation}
P_0(\xi) = \eta \mathrm{e}^{-\eta \xi},
\end{equation}
with
\begin{equation}
{\eta}=\frac{2\omega\E[B_2']^2}{2\frac{\lambda_1}{\mu_1^2}+\frac{1}{\E[B_2']}\left(1-\frac{\lambda_1}{\mu_1}\right)\E[B_2'^2]}.
\label{eta}
\end{equation}
\end{theorem}

A natural choice for $\omega$ is $\omega = 1/\E[B_2']$, leading to the commonly used scaling $\lim_{\delta\downarrow0}\delta N_2=\lim_{u_2\uparrow1}(1-u_2)N_2$. We refer to Section \ref{critval} for a discussion on this choice and alternatives. By applying the multiclass distributional law of Bertsimas and Mourtzinou \cite{bertsimasmourtzinou} it directly follows that the scaled waiting time at $Q_2$ also follows an exponential distribution with parameter $\lambda_2 \eta$.

Finally, the above has the following immediate consequence for the joint (scaled) queue-length distribution in HT.

\begin{theorem}
For $\lambda_1/k_1<\lambda_2/k_2$ and $\lambda_2=(1-\rho_1-\delta)/{\E[B_2']}$, we have:
\begin{equation}
\lim_{\delta\downarrow0}\P[N_1\leq n_1, \delta N_2\leq \xi] = \mathcal{L}(n_1)\left(1-\mathrm{e}^{-\eta \xi}\right),\label{mainresult}
\end{equation}
where $\eta$ is given by \eqref{eta}, and $\mathcal{L}(\cdot)$ is the cumulative probability distribution of the steady-state queue length of a single-server multiple-vacation queue with parameters $\lambda_1$ and $\mu_1$, with $k_1$-limited service. The vacation distribution in this system equals the convolutions of $k_2$ service-time distributions in $Q_2$ and of the switch-over times $S_1$ and $S_2$.
\end{theorem}

\section{Discussion}  \label{discussions}

\subsection{Interpretation}
The main result \eqref{mainresult} has the following intuitively appealing interpretation, which we derive heuristically below:
\begin{compactenum}
\item The number of customers in the stable queue has the same distribution as the number of customers in a $k$-limited vacation system.\label{propone}
\item The scaled number of customers in the critically loaded queue is exponentially distributed with parameter $\eta$.\label{proptwo}
\item The number of customers in the stable queue and the (scaled) number of customers in the critically loaded queue are independent.\label{propthree}
\end{compactenum}

Property \ref{propone} follows from the fact that if $Q_2$ is in HT, then exactly $k_2$ customers are served at this queue during each
cycle.

For Property \ref{proptwo}, we take a look at an alternative $k$-limited polling system with the same arrival rates, \textit{no} switch-over times, the original service times  at $Q_1$ and adjusted service times at $Q_2$ with mean $\mathbb{E}[B_2']$ and variance $\V[B_2']$. Since in the original system  the two switch-over periods $S_1$ and $S_2$ are always accompanied with precisely $k_2$ service times at $Q_2$ in HT, the first two moments of the  amount of time the server is utilized - either due to serving or switching - in the original and alternative system is obviously identical. In \cite{boonperturbation} it is proven that the distribution of the scaled workload within the alternative system without switch-over times  - and, thus, the original system -  equals the scaled amount of work in an $M/G/1$ queue in which the two customer classes are combined into one customer
class. Based on standard heavy-traffic results for the $M/G/1$ queue this implies that the
distribution of the scaled total workload converges to an exponential distribution with mean $ \rho   \mathbb{E}[B^2] / 2 \mathbb{E}[B]$, where $\mathbb{E}[B] =(\rho_1 + \lambda_2\mathbb{E}[B_2'])/ (\lambda_1+\lambda_2)$  and  $\mathbb{E}[B^2] = (2 \lambda_1/ \mu_1^2 + \lambda_2 (\V[B_2'] + \mathbb{E}[B_2']^2) ) / (\lambda_1+\lambda_2)$.  Since in HT almost all customers are located in $Q_2$ and $\lambda_2 \uparrow (1-\rho_1)/\mathbb{E}[B_2']$, the scaled number of customers in $Q_2$ is exponentially distributed with parameter $\eta$.

Property \ref{propthree} is due to the the time-scale separation in HT, i.e., the dynamics of the stable queue evolve at a much faster time scale than the dynamics of the critically loaded queue.

\subsection{Critical value and choice of \boldmath $\omega$}\label{critval}
In this note, we take the limit by increasing the arrival rate of $Q_2$ until it reaches its critical value while keeping the other parameters fixed. However, one could consider more general ways of varying the arrival rates in order to let $Q_2$ become critically loaded. To this end, we introduce $\lambda_1^*$ and $\lambda_2^*$ such that
$\lambda_1^*/\mu_1+\lambda_2^* \E[B_2']= 1$. Additionally, we assume that
\begin{equation}
\frac{\lambda_1^*}{k_1}<\frac{1}{\frac{k_1}{\mu_1}+k_2 \E[B_2']},
\label{stability2alt}
\end{equation}
or equivalently:
\begin{equation}
\frac{\lambda_2^*}{k_2}>\frac{1}{\frac{k_1}{\mu_1}+ k_2 \E[B_2']}.
\end{equation}
We now let $\lambda_1\rightarrow\lambda_1^*$ and $\lambda_2\rightarrow\lambda_2^*$ for $\delta\downarrow 0$, with
\begin{equation}
\frac{\lambda_1}{\mu_1} + \lambda_2 \E[B_2'] = 1-\delta\omega^*,\qquad \omega^* > 0, 0<\delta \ll 1.\label{perturbationalt}
\end{equation}
Any arbitrary way in which we let $\lambda_1$ and $\lambda_2$ approach respectively $\lambda_1^*$ and $\lambda_2^*$, for $\delta\downarrow 0$, will cause $Q_2$ to become critically loaded (because of assumption \eqref{stability2alt}). All results obtained in this paper will still be valid, by choosing $\omega^* = \omega \E[B_2']$. 

In Section \ref{mainresults} we argued that $\omega=1/\E[B_2']$ is a good choice, because it leads to $\delta N_2 = (1-u_2)N_2$. To see this, we first introduce the notation $\lambda_2^\textit{crit}=(1-\rho_1)/\E[B_2']$ and rewrite:
\begin{equation}
\lim_{\delta\downarrow0}\delta N_2=\lim_{\lambda_2\uparrow \lambda_2^\textit{crit}}\left(\frac{1}{\omega}(\lambda_2^\textit{crit}-\lambda_2)\right)N_2.
\label{scalinglambda2}
\end{equation}
Substitution of $\omega=1/\E[B_2']$ yields
\begin{align*}
\lim_{\lambda_2\uparrow \lambda_2^\textit{crit}} \left(\frac{1}{\omega}(\lambda_2^\textit{crit}-\lambda_2)\right)N_2
& = \lim_{\lambda_2\uparrow \lambda_2^\textit{crit}}\left((1-\rho_1)-\lambda_2\E[B_2'])\right)N_2\\
&= \lim_{u_2\uparrow 1}\left(1-u_2\right)N_2,
\end{align*}
which has an appealing form. From \eqref{scalinglambda2} one can immediately see that $\omega=1$ results in taking the limit of $(\lambda_2^\textit{crit}-\lambda_2)N_2$, which is also quite commonly used (see for example \cite{tedijanto}).

Finally, we consider increasing $\lambda_1$ and $\lambda_2$ simultaneously while keeping their \emph{ratio} fixed. Denote by $\Lambda=\lambda_1+\lambda_2$ the total arrival rate, and let $\lambda_i=\hat\lambda_i\Lambda$ $(i=1,2)$. As before, without loss of generality, we assume that $Q_2$ becomes unstable when increasing $\Lambda$. Taking $\omega=\hat\lambda_2$ yields
\begin{align*}
\lim_{\lambda_2\uparrow \lambda_2^\textit{crit}}\left(\frac{1}{\omega}(\lambda_2^\textit{crit}-\lambda_2)\right)N_2
&=\lim_{\Lambda\uparrow \Lambda^\textit{crit}}(\Lambda^\textit{crit}-\Lambda)N_2,
\end{align*}
where $\Lambda^\textit{crit}$ is the total arrival rate at which the second queue becomes critically loaded.

\subsection{Two critically loaded queues}
In the current note we have analysed the heavy-traffic behaviour under the condition $\lambda_1/k_1<\lambda_2/k_2$, i.e., only $Q_2$ becomes critically loaded. The limiting regime, in which both queues become saturated simultaneously $(\lambda_1/k_1=\lambda_2/k_2)$, shows fundamentally different system behaviour. In this case, it is not sufficient anymore to use a scaling that implies that exactly $k_2$ customers are served at $Q_2$ during each cycle, i.e., the probability that there are served less than $k_2$ customers cannot be neglected, when analyzing the asymptotic behaviour of $Q_1$.

\section{Extensions}

The presented framework can be generalized in several areas without fundamentally complicating the analysis.

\subsection{General distributions}  Using the insights in the heavy-traffic system behavior obtained from the perturbation analysis, a HT limit theorem for systems with general renewal arrival processes and general service and switch-over times can be conjectured. The interchange between service and switch-over times proves its value again, since it makes it possible to link the workload of the system to the workload of a $GI/G/1$ queue, for which the HT behavior is well-studied. Subsequently, the limiting distribution of the (scaled) number of customers in $Q_2$ can be derived. Let $A_i$, $B_i$, and $S_i$ denote the generic interarrival times, service times, and switch-over times for $i=1,2$, which can be generally distributed. Furthermore, let
\begin{equation*}
\E[B_2']=\E[B_2]+\frac{\E[S]}{k_2},\quad
\V[B_2']=\V[B_2]+\frac{\V[S]}{k_2},\label{meanvarB2general}
\end{equation*}
which is a generalisation to general service times of \eqref{meanvarB2}, and let
\[
\lambda_2^\textit{crit} = (1-\rho_1)/\mathbb{E}[B_2']
\]
denote the limiting arrival intensity of customers at $Q_2$. We scale the interarrival-time distribution $A_2$ by changing its mean, while keeping the coefficient of variation fixed. We can formulate the following conjecture.

\begin{conjecture}\label{conjecturegeneraldistributions}
The HT limit of the scaled queue length $(1-u_2)N_2$ as $u_2\uparrow 1$ is exponentially distributed with mean
\begin{align}
\frac{1}{\eta} =& \frac{1}{2\E[B_2']}\Big(\lambda_1(\V[B_1]+\SCV{A_1}\E[B_1]^2)\nonumber\\
&+\lambda_2^\textit{crit}(\V[B_2']+\SCV{A_2}\E[B_2']^2)]\Big),\label{etageneral}
\end{align}
where $\SCV{A_i}$ denotes the squared coefficient of variation of $A_i$.
\end{conjecture}

\subsection{\boldmath$N$ ($>2$) queue polling system }\label{morethannqueues}
One could readily extend the analysis and results of the present note to an $N$($>2$) queue polling system with one critically loaded queue. In this case the stable queues have the same joint queue-length distribution as in an $N-1$ queue polling model with an extended switch-over time, while the scaled queue-length distribution of the critically loaded queue follows again an exponential distribution.

\section{Numerical results}

This section presents four numerical examples illustrating the validity and applicability of our proven limiting results and key concepts:
\begin{enumerate}
\item (Joint) queue length distributions;
\item Variance absorption result;
\item Extension to general distributions;
\item Extension to $N$($>2$) queue polling system.
\end{enumerate}

\noindent\textbf{Example 1.} In the first example we consider three polling models, referred to as $P1$, $P2$, and $P3$. In each model we have $\lambda_1=0.06$, $\mu_1=\mu_2=0.5$, $\sigma_1=\sigma_2=0.5$, $k_1=2$ and $k_2=3$. The only difference between the three models are the values of $\lambda_2$, which are  $0.25$, $0.258$, and $0.263$, respectively. The corresponding values for $u_2$ are $0.953$, $0.98$, and $0.997$. Table \ref{tbl:ql1} shows the simulated queue-length probabilities for $Q_1$ for these three models, and for the corresponding vacation model (model $V$). It can clearly be seen that the queue-length probabilities converge to those of the vacation system as stated in Theorem \ref{theorem1}.

\begin{table}[b]
\[
\begin{array}{|c|cccccc|}
\hline
\text{Model} & p_0 & p_1 & p_2 & p_3 & p_4 & p_5 \\
\hline
P1 &        0.589 & 0.275 & 0.095 & 0.029 & 0.009 & 0.003 \\
P2 &        0.585 & 0.277 & 0.096 & 0.030 & 0.009 & 0.003 \\
P3 &        0.583 & 0.278 & 0.097 & 0.030 & 0.009 & 0.003 \\
\hline
V &   0.582 & 0.278 & 0.097 & 0.030 & 0.009 & 0.003 \\
\hline
\end{array}
\]
\caption{\boldmath Simulated values for $p_k := \P[N_1=k]$.}
\label{tbl:ql1}
\end{table}

Figure \ref{fig:ql2} shows how in the simulation the scaled queue-length $(1-u_2)N_2$ in model $P1$ closely follows the exponential distribution with parameter $\eta=1.289$ as expected from Theorem \ref{theorem2}. For models $P2$ and $P3$ the differences between the simulation and the closed-form asymptotics are hardly noticeable, which is the reason why we have omitted them in the figure.

\begin{figure}[t]
\begin{center}
\includegraphics[height=4cm]{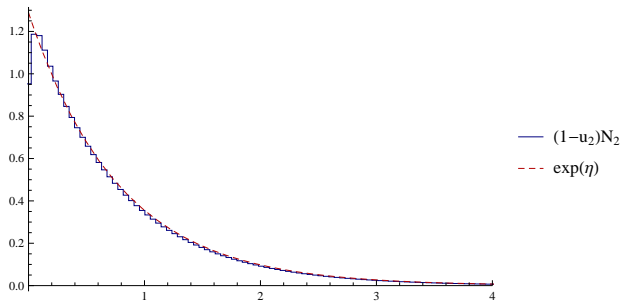}
\end{center}
\caption{\boldmath The densities of $(1-u_2)N_2$ of model $P1$ and an exponential distribution with parameter $\eta$.}
\label{fig:ql2}
\end{figure}

Finally, we have calculated the  correlation coefficient between the simulated number of customers in $Q_1$ and the scaled queue length of $Q_2$. For $P1$, $P2$, and $P3$ the values are  0.048, 0.022, and 0.004, respectively. These values converge to $0$, which follows from Property \ref{propthree}, implying that the amount of ``memory'' of the stable queue asymptotically vanishes compared to that of the critically loaded queue. \\

\noindent\textbf{Example 2.} As a by-product of our analysis we have observed that an interchange exists among the first two moments in service and switch-over times in systems with and without switch-over times such that the HT limits (for $Q_2$) are identical in both systems. That is, by noting that in HT every cycle precisely $k_2$ are served,  the second terms in $ \E[B_2']$ and $ \V[B_2']$  apportion the first two moments of the switch-over times to each of the customers fairly. It is important to note that this interchange is just one of the possibilities, i.e., a variety of trade-offs exists between service and switch-over times, all of which yield the same asymptotic distributions, as long as the terms \eqref{meanvarB2} remain unchanged.

In this second example, we want to illustrate this generalization of the variance absorption result for branching systems \cite{cooper} with a numerical case. We start by considering again the three polling models $P1$, $P2$, and $P3$, but now we construct alternative models where (a portion of) the switch-over times is transferred to the service times of the customers in $Q_2$. In more detail, we construct corresponding models with modified service times and switch-over times, satisfying
\begin{align*}
\E[B_2]& =\frac{1}{\mu_2}+\alpha\frac{1/\sigma_1+1/\sigma_2}{k_2}, \\
\V[B_2]& =\frac{1}{\mu_2^2}+\alpha\frac{1/\sigma_1^2+1/\sigma_2^2}{k_2},\\
\E[S_1]& = (1-\alpha)/\sigma_1, \qquad
\V[S_1]=(1-\alpha)/\sigma_1^2, \\
\E[S_2]& = (1-\alpha)/\sigma_2, \qquad
\V[S_2]=(1-\alpha)/\sigma_2^2,
\end{align*}
choosing $\alpha=1/2$ for models $R1$, $R2$, $R3$, and $\alpha=1$ for models $S1$, $S2$, $S3$. Note that these latter three models are polling models \emph{without} switch-over times and that $u_2=\rho$ in this case. The means, standard deviations, and squared coefficients of the scaled queue length $(1-u_2)N_2$ are obtained using simulation, and given in Table \ref{tbl:ql2example2}. These values clearly confirm that the asymptotic distributions are identical, independent of the chosen interchange, and that they converge to an exponential distribution with mean $1/\eta=1/1.289=0.776$.
If the aforementioned interchange keeps the \textit{distribution} of the sum of $k_2$ service times in $Q_2$ and the switch-over times unaffected, then the HT limit of $Q_1$ remains unchanged as well.
\begin{table}[t]
\[
\begin{array}{|l|rrr|}
\hline
(1-u_2)N_2 & R1 & R2 & R3 \\
\hline
\text{Mean}    &  0.742 & 0.762 & 0.778  \\
\text{St.dev.} &  0.757 & 0.766 & 0.781  \\
\text{CV}      &  1.020 & 1.006 & 1.005  \\
\hline
\hline
(1-\rho)N_2 & S1 & S2 & S3  \\
\hline
\text{Mean}     &  0.730 & 0.757 & 0.770 \\
\text{St.dev.}  &  0.756 & 0.768 & 0.768 \\
\text{CV}       &  1.035 & 1.015 & 0.996 \\
\hline
\end{array}
\]
\caption{Simulated means, standard deviations, and coefficients of variation of\boldmath $(1-u_2)N_2$ for Example 2.}
\label{tbl:ql2example2}
\end{table}
\\
\
\\
\noindent\textbf{Example 3.}
In this example we use simulation to validate Conjecture \ref{conjecturegeneraldistributions} on the HT asymptotics for systems with general renewal arrival processes and general service and switch-over times. In many communication systems (some of these) distributions may be heavy-tailed, stressing the importance of knowing the limiting queue-lengths for general distributions. We take the same setting as in Example 1, but with different distributions. Table \ref{tbl:example3input} gives an overview of the various interarrival-time, service-time, and switch-over time distributions. Note that their means are equal to those in Example 1, where everything was assumed to be exponentially distributed. 
Also note that we have added a fourth system, P4, with $\lambda_2=0.2.635$, for reasons that become clear later.

\begin{table}[!htb]
\begin{center}
\begin{tabular}{|l|c|c|}
\hline
Distribution & Queue 1 & Queue 2\\
\hline
Interarrival times & Gamma$(0.5, 0.03)$ & Pareto$(2/(3\lambda_2), 3)$ \\
\hline
Service times & Uniform$(0, 4)$ & Pareto$(1.17, 2.4)$ \\
\hline
Switch-over times & Constant$(2)$ & Constant$(2)$ \\
\hline
\end{tabular}
\caption{Settings for Example 3.}
\label{tbl:example3input}
\end{center}
\end{table}

We have simulated the queue-length distributions of each of these four systems and compare them to the limiting values, based on the exponential distribution that arises in the HT limit. The parameter of this exponential distribution, which can be computed using \eqref{etageneral}, is $\eta\approx2.527$. Table \ref{tbl:ql2example3} gives the mean, standard deviations, and coefficients of variation of the scaled queue lengths for each of the models. From this table we conclude that the first two moments indeed approach those of an exponential distribution, albeit in a slower pace than Example 1, which had only exponential distributions. This slower speed of convergence is the reason for including a fourth model, P4, with an even higher utilization. Figures \ref{fig:ql3}(a)--(b) shows the simulated scaled queue length probabilities for P1, \dots, P4, and the density of the limiting distributions, again confirming the convergence to the exponential distribution.

\begin{table}[t]
\[
\begin{array}{|l|rrrr|}
\hline
(1-u_2)N_2 & P1 & P2 & P3 & P4\\
\hline
\text{Mean}    &    7.5 & 18.4 & 116.3 & 232.5 \\
\text{St.dev.} &    30.5 & 38.6 & 192.1 & 281.8 \\
\text{CV}      &    4.0 & 2.1 & 1.7 & 1.2 \\
\hline
\end{array}
\]
\caption{Simulated means, standard deviations, and coefficients of variation for Example 3.}
\label{tbl:ql2example3}
\end{table}

\begin{figure}[!htb]
\begin{tabular}{cc}
\includegraphics[width=0.45\linewidth]{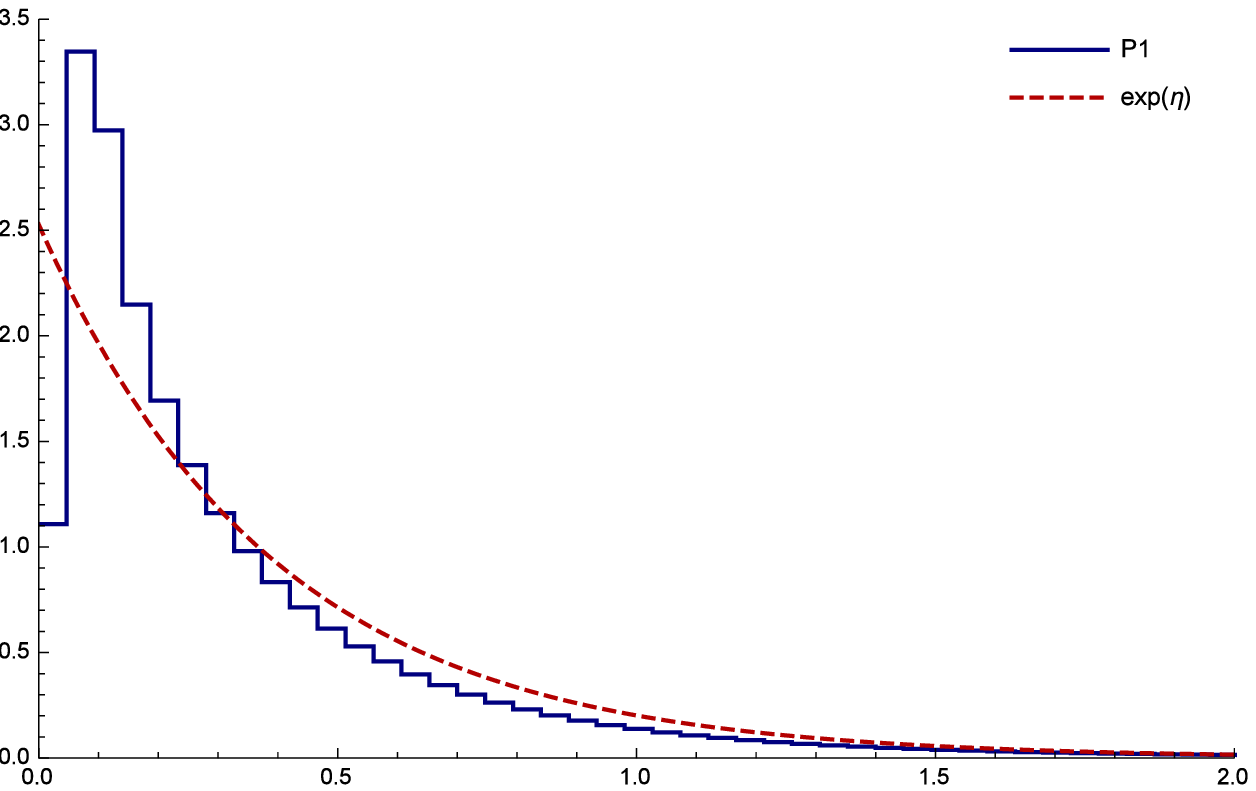} & \includegraphics[width=0.45\linewidth]{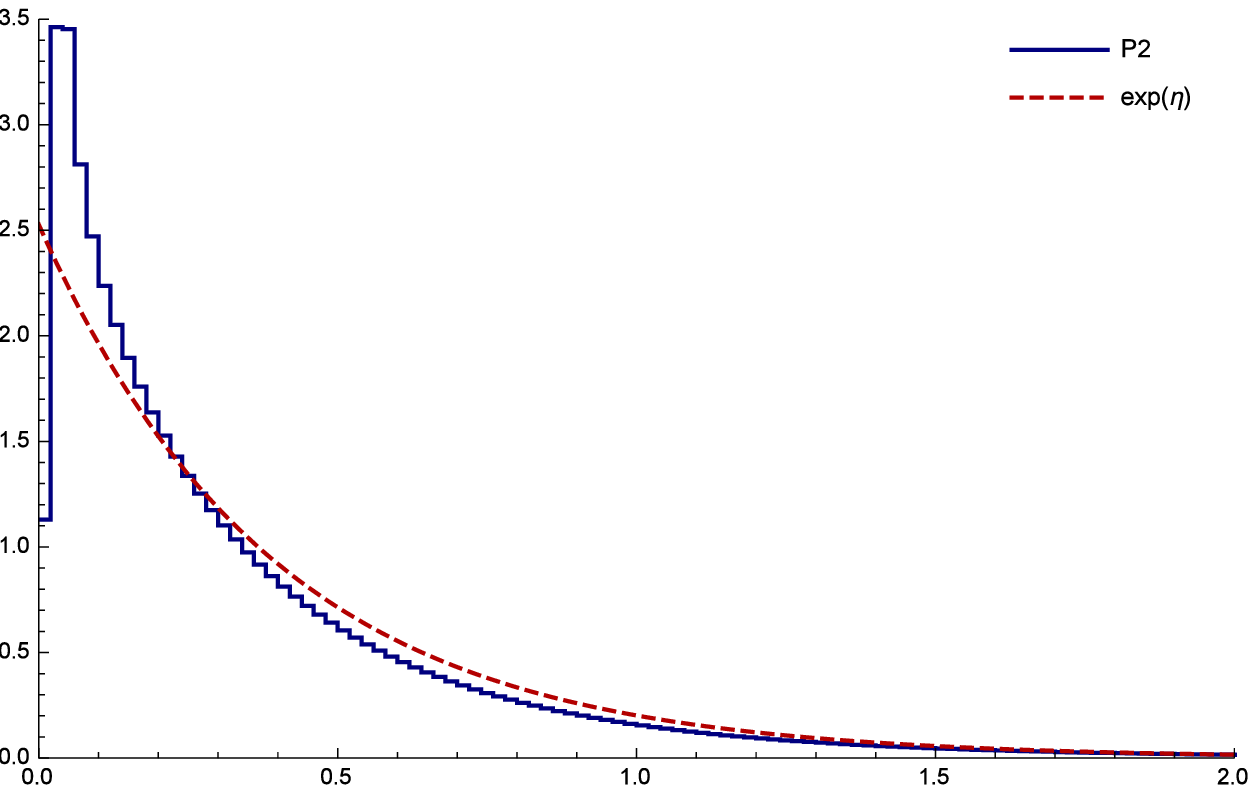}  \\
(a) P1 $(u_2=0.953)$& (b) P2 $(u_2=0.980)$\\
\includegraphics[width=0.45\linewidth]{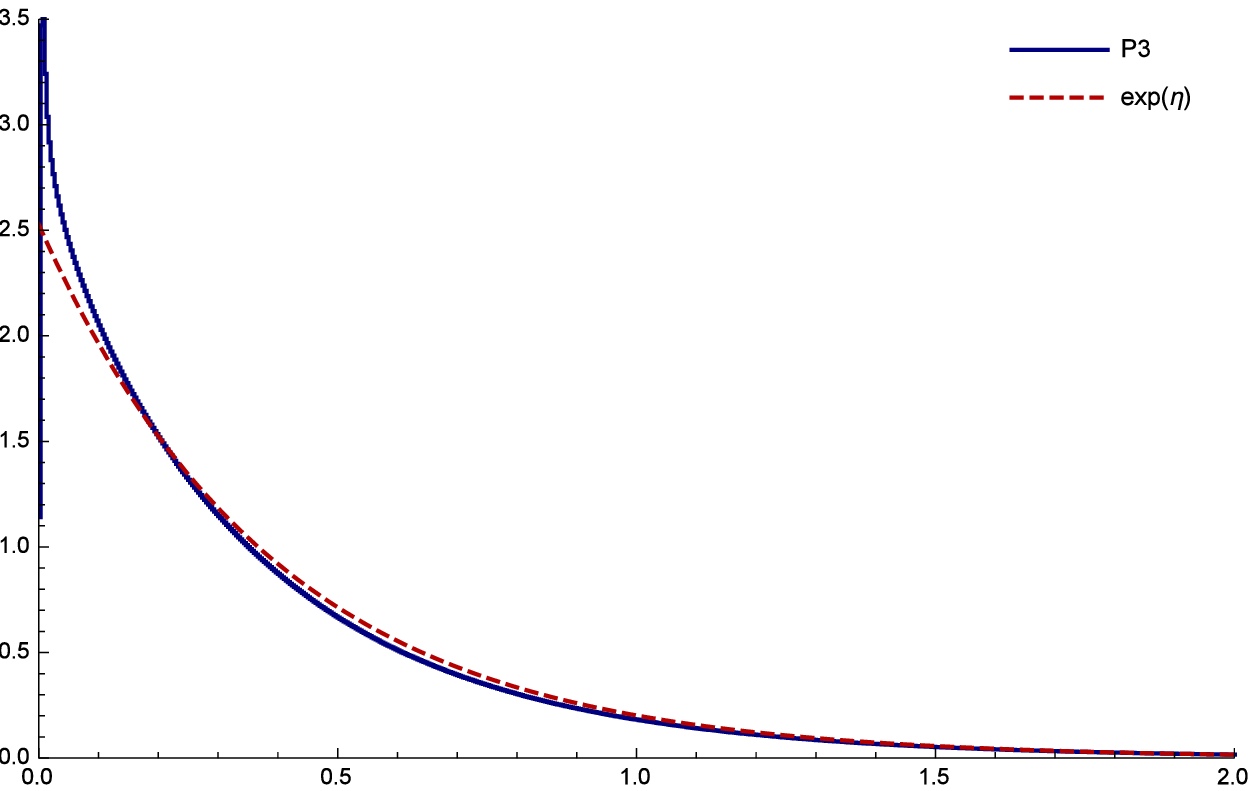} & \includegraphics[width=0.45\linewidth]{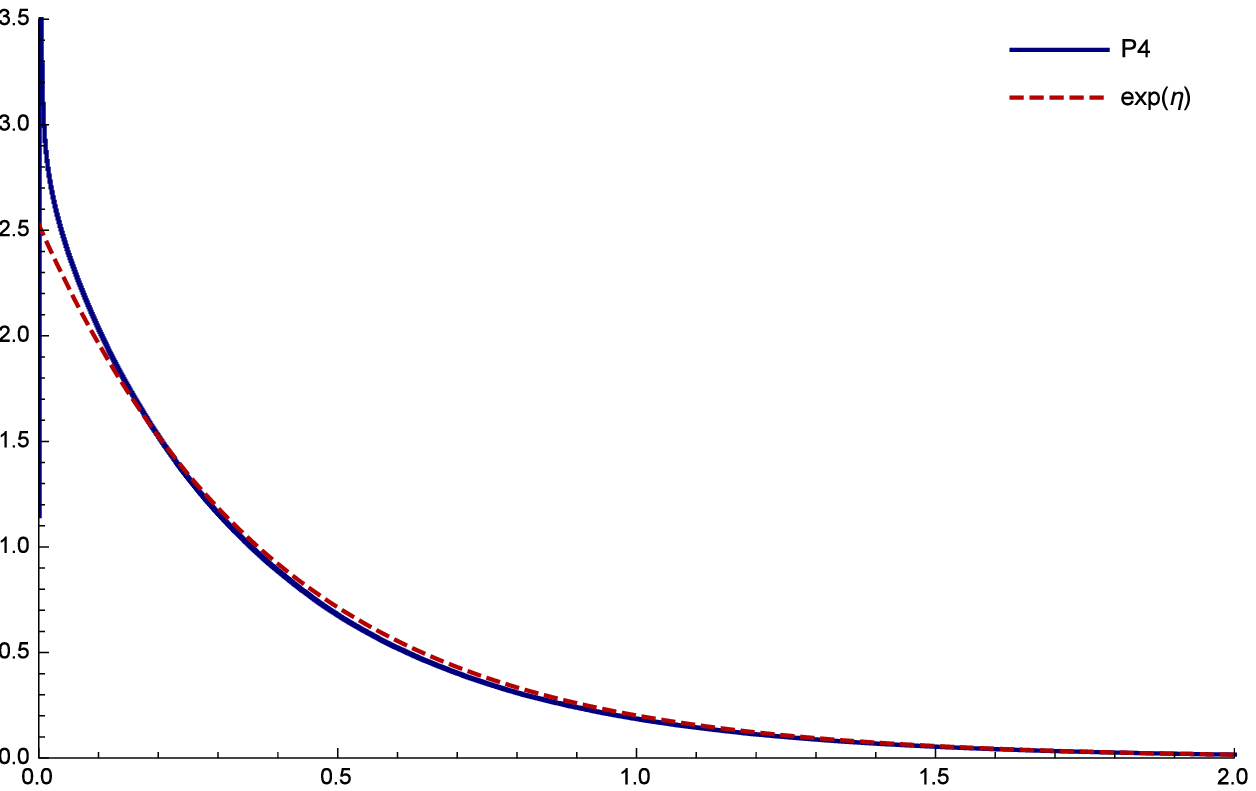}\\
(c) P3 $(u_2=0.997)$& (d) P4 $(u_2=0.998)$\\
\end{tabular}
\caption{\boldmath Simulated and theoretical scaled queue lengths for Example 3.}
\label{fig:ql3}
\end{figure}

\
\\
\noindent\textbf{Example 4.}
In this final example we illustrate how queue lengths behave in a polling system with $N (> 2)$ queues, each with $k$-limited service, when \emph{all} arrival rates are increased simultaneously, as discussed in Section \ref{critval}. We consider a system with $N=4$ queues and exponentially distributed service times, interarrival times, and switch-over times with $\mu_i = 1/2$, $\sigma_i=1/3$, for $i=1,2,3,4$. The service limits are $k_1=7$, $k_2=k_3=6$, and $k_4=5$.
Let $\Lambda=\lambda_1+\dots+\lambda_4$ denote the total arrival rate. We will increase $\Lambda$ until, one by one, \emph{all} of the queues become unstable. While increasing $\Lambda$, we keep the ratio between the four arrival rates fixed. Let $\hat\lambda_i$ denote the proportion of the total arrivals that is routed to $Q_i$, i.e.
\begin{equation}
\lambda_i=\hat\lambda_i\Lambda,\label{lambdai}
\end{equation}
with $\hat\lambda_i=i/10$ in this example. The  impact of increasing the total arrival rate on the system behavior is graphically illustrated in Figure \ref{fig:ql4}, which depicts simulated mean queue lengths $\E[N_i(\Lambda)]$ as a function of the total arrival rate $\Lambda$.
\begin{figure}[!hbtp]
\begin{center}
\begin{tikzpicture}
\node[anchor=south west,inner sep=0] at (0,0) {\includegraphics[width=0.9\linewidth]{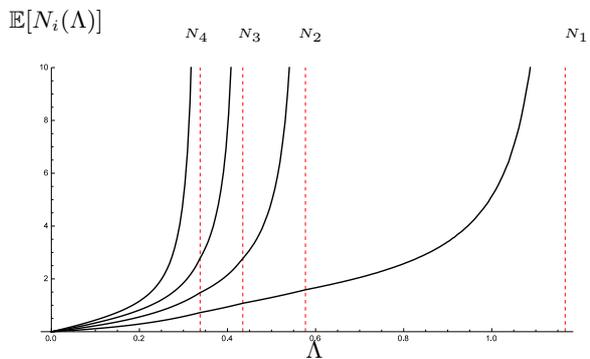}};
\node[anchor=south west,inner sep=0] at (-0.4,4.1) {\small $\E[N_i(\Lambda)]$};
\node[anchor=south west,inner sep=0] at (1.9,4.0) {\tiny $N_4$};
\node[anchor=south west,inner sep=0] at (2.6,4.0) {\tiny $N_3$};
\node[anchor=south west,inner sep=0] at (3.4,4.0) {\tiny $N_2$};
\node[anchor=south west,inner sep=0] at (6.9,4.0) {\tiny $N_1$ };
\node[anchor=south west,inner sep=0] at (3.5,-0.2) {\small $\Lambda$ };
\end{tikzpicture}
\caption{The mean queue lengths in Example 4 as functions of the total arrival rate $\Lambda$. The four vertical dashed lines correspond to respectively $\Lambda_4^\textit{crit}, \Lambda_3^\textit{crit}, \Lambda_2^\textit{crit}, \Lambda_1^\textit{crit}$ (from left to right).}
\label{fig:ql4}
\end{center}
\end{figure}

\noindent The following interesting observations can be made:
\begin{itemize}
\item Obviously, for $\Lambda = 0$ all mean queue lengths are 0.
\item As $\Lambda$ is increased, all mean queue lengths increase. Queue 1 has the lowest utilization, causing the mean queue length to increase slowest. Queue 4 has the highest utilization and, as a consequence, is the first queue to become unstable, which happens at $\Lambda_4^\textit{crit}=25/74\approx0.34$ (the leftmost dashed vertical line). This value can be computed via \eqref{eqn:ui} and \eqref{lambdai}.
\item Note that, although $Q_4$ is unstable at $\Lambda=\Lambda_4^\textit{crit}$, the other queues are stable. We can increase $\Lambda$ further, causing the mean queue lengths of the other queues to increase, while $Q_4$ remains unstable.
\item When $\Lambda$ approaches the value of $\Lambda_3^\textit{crit}=10/23\approx0.43$, $Q_3$ becomes unstable as well. Now $Q_1$ and $Q_2$ are the only stable queues.
\item As we continue to increase the total arrival rate, $Q_2$ and, eventually, also $Q_1$ become unstable. This happens at respectively $\Lambda_2^\textit{crit}=15/26\approx0.58$ and $\Lambda_1^\textit{crit}=7/6\approx1.17$.
\item The mean queue lengths (as functions of $\Lambda$) are not differentiable at $\Lambda=\{\Lambda_1^\textit{crit}, \dots, \Lambda_4^\textit{crit}\}$.
\end{itemize}
This particular behavior can be explained using the results from this paper. As discussed in Section \ref{morethannqueues}, from the viewpoint of customers in queues 1, 2, and 3, the four-queue polling system behaves like a three-queue polling system when $Q_4$ is unstable, i.e. for $\Lambda_4^\textit{crit}<\Lambda<\Lambda_3^\textit{crit}$. In this so-called ``corresponding system'' the switch-over time between $Q_3$ and $Q_1$ is the sum of the switch-over times $S_3+S_4$ in the original model, plus $k_4$ service times that are $\exp(\mu_4)$ distributed. In this three-queue polling model, the next queue to become unstable is $Q_3$, at $\Lambda_3^\textit{crit}=10/23\approx0.43$, making the system behave like a two-queue polling model with an even longer switch-over time between $Q_2$ and $Q_1$. This continues, with $Q_2$ becoming unstable at $\Lambda=15/26\approx0.58$ and, finally, $Q_1$ being the last queue to become unstable when the total arrival rate reaches the value of $\Lambda_1^\textit{crit}=7/6\approx1.17$.

The simulation results can be used to verify our claims that, near the point where one of the queues becomes unstable, (a) the number of customers in the stable queues are distributed as in a system with that particular queue removed and replaced by an extra long switch-over time, and (b) the scaled queue length of the unstable queue is exponentially distributed. Therefore, we consider the original polling system near each of the critical values, namely at $\Lambda=0.995\Lambda_i^\textit{crit}$, for all $i$, and compare the simulation results to the simulation results of the ``corresponding'' polling systems. Table \ref{tbl:ql4}(a) shows the simulated queue length probabilities of the stable queues. For $\Lambda \approx \Lambda_4^\textit{crit}$ queue 4 becomes unstable and the system starts behaving as a three-queue polling system with an extra long switch-over time $S_3$. The simulated queue length probabilities of this ``corresponding system'' are listed in Table \ref{tbl:ql4}(b). It is readily seen that the results are completely equivalent, the extremely minor differences being attributed to simulation inaccuracies. This has been repeated near $\Lambda \approx \Lambda_3^\textit{crit}$
and $\Lambda \approx \Lambda_2^\textit{crit}$,
and the same conclusion can be drawn.
\begin{table}[!htb]
\[
\begin{array}{|l|cccccc|}
\hline
& \multicolumn{6}{|c|}{\Lambda \approx \Lambda_4^\textit{crit}}\\
\hline
\text{Queue} & p_0 & p_1 & p_2 & p_3 & p_4 & p_5  \\
\hline
Q_1 & 0.534 & 0.292 & 0.118 & 0.040 & 0.012 & 0.003 \\
Q_2 & 0.325 & 0.281 & 0.186 & 0.106 & 0.055 & 0.026 \\
Q_3 & 0.186 & 0.205 & 0.177 & 0.136 & 0.098 & 0.067 \\
\hline
\multicolumn{7}{}{}\\
\hline
& \multicolumn{6}{|c|}{\Lambda \approx \Lambda_3^\textit{crit}}\\
\hline
\text{Queue} & p_0 & p_1 & p_2 & p_3 & p_4 & p_5  \\
\hline
Q_1 & 0.410 & 0.301 & 0.166 & 0.076 & 0.031 & 0.011 \\
Q_2 & 0.177 & 0.200 & 0.179 & 0.142 & 0.103 & 0.071 \\
\hline
\multicolumn{7}{}{}\\
\hline
& \multicolumn{6}{|c|}{\Lambda \approx \Lambda_2^\textit{crit}}\\
\hline
\text{Queue} & p_0 & p_1 & p_2 & p_3 & p_4 & p_5  \\
\hline
Q_1 & 0.296 & 0.273 & 0.195 & 0.119 & 0.064 & 0.031 \\
\hline
\end{array}
\]
\centering
(a) Original polling system\\

\[
\begin{array}{|l|cccccc|}
\hline
& \multicolumn{6}{|c|}{\Lambda \approx \Lambda_4^\textit{crit}}\\
\hline
\text{Queue} & p_0 & p_1 & p_2 & p_3 & p_4 & p_5  \\
\hline
Q_1 & 0.534 & 0.292 & 0.118 & 0.040 & 0.012 & 0.003 \\
Q_2 & 0.326 & 0.281 & 0.185 & 0.106 & 0.055 & 0.026 \\
Q_3 & 0.187 & 0.206 & 0.176 & 0.135 & 0.097 & 0.067 \\
\hline
\multicolumn{7}{}{}\\
\hline
& \multicolumn{6}{|c|}{\Lambda \approx \Lambda_3^\textit{crit}}\\
\hline
\text{Queue} & p_0 & p_1 & p_2 & p_3 & p_4 & p_5  \\
\hline
Q_1 & 0.411 & 0.301 & 0.165 & 0.076 & 0.030 & 0.011 \\
Q_2 & 0.177 & 0.201 & 0.179 & 0.142 & 0.103 & 0.070 \\
\hline
\multicolumn{7}{}{}\\
\hline
& \multicolumn{6}{|c|}{\Lambda \approx \Lambda_2^\textit{crit}}\\
\hline
\text{Queue} & p_0 & p_1 & p_2 & p_3 & p_4 & p_5  \\
\hline
Q_1 & 0.297 & 0.273 & 0.195 & 0.119 & 0.064 & 0.031 \\
\hline
\end{array}
\]
\centering
(b) Corresponding polling systems\\

\caption{\boldmath Marginal queue lengths of the stable queues in the original polling and the corresponding systems. The simulated probabilities $\P(N_i=m)$, denoted by $p_m$, are given for $m=0, 1, 2, \dots, 5$.}\vspace*{-4ex}
\label{tbl:ql4}
\end{table}

Verifying that the scaled queue lengths of the unstable queues are exponentially distributed, which is our second claim, requires determining the parameters of these distributions  as a first step. These parameters can be computed similar to \eqref{etageneral}, but with $B_2'$ modified to include the service times and switch-over times of \emph{all} overloaded queues, and $B_1$ the service time of an arbitrary customer in \emph{any} of the stable queues. The values of $\eta_i$, $i=1,2,3,4$, equal $1.364$, $2.875$, $3.739$ and $4.529$ taking $\omega=\hat\lambda_i$. The latter is the appropriate value when considering the scaling $(\Lambda_i^\textit{crit}-\Lambda)N_i$ as suggested in Section \ref{critval}. Figures \ref{ql4b}(a)--(d) depict the simulated probability density functions of the scaled queue lengths $(\Lambda_i^\textit{crit}-\Lambda)N_i$, for $i=1,\dots,4$ at $\Lambda=0.995\Lambda_i^\textit{crit}$, i.e. very close to instability. It is clear that the simulated values are extremely close to the theoretical values.


\begin{figure}[!htb]
\begin{tabular}{cc}
\includegraphics[width=0.45\linewidth]{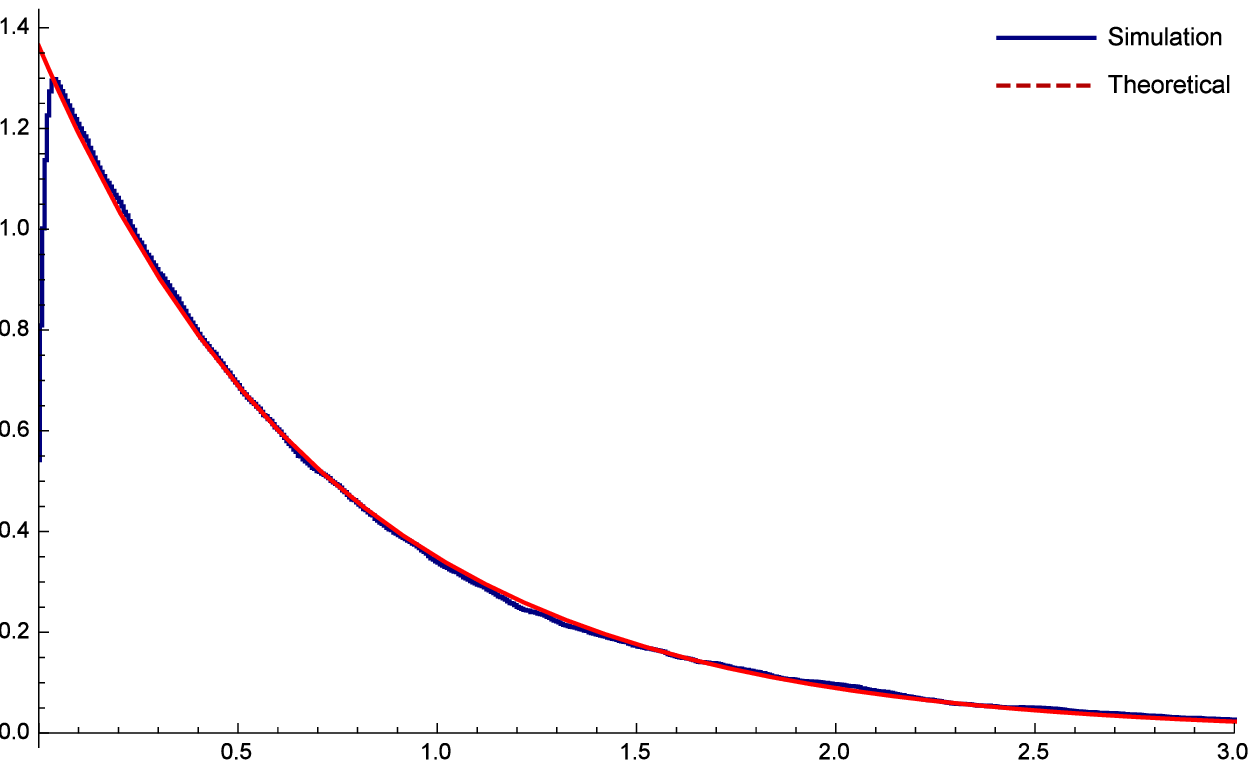} & \includegraphics[width=0.45\linewidth]{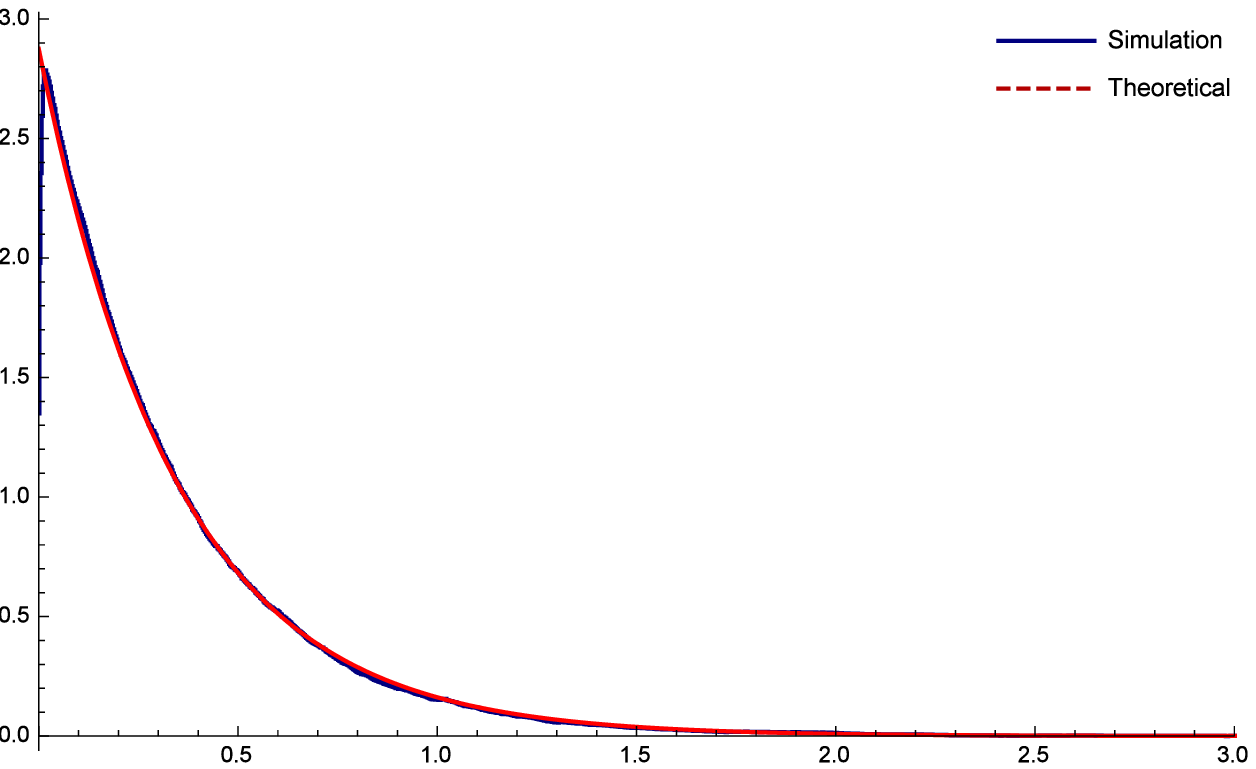}  \\
(a) Queue 1 & (b) Queue 2 \\
\includegraphics[width=0.45\linewidth]{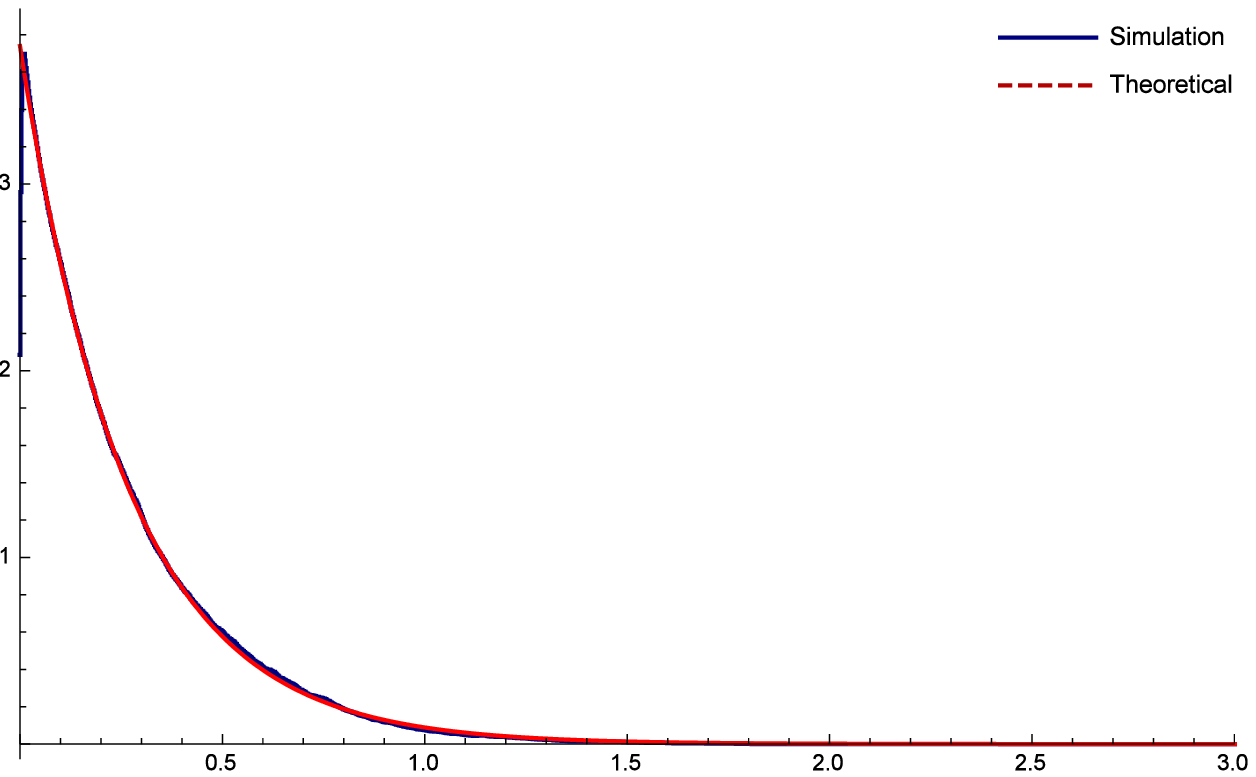} & \includegraphics[width=0.45\linewidth]{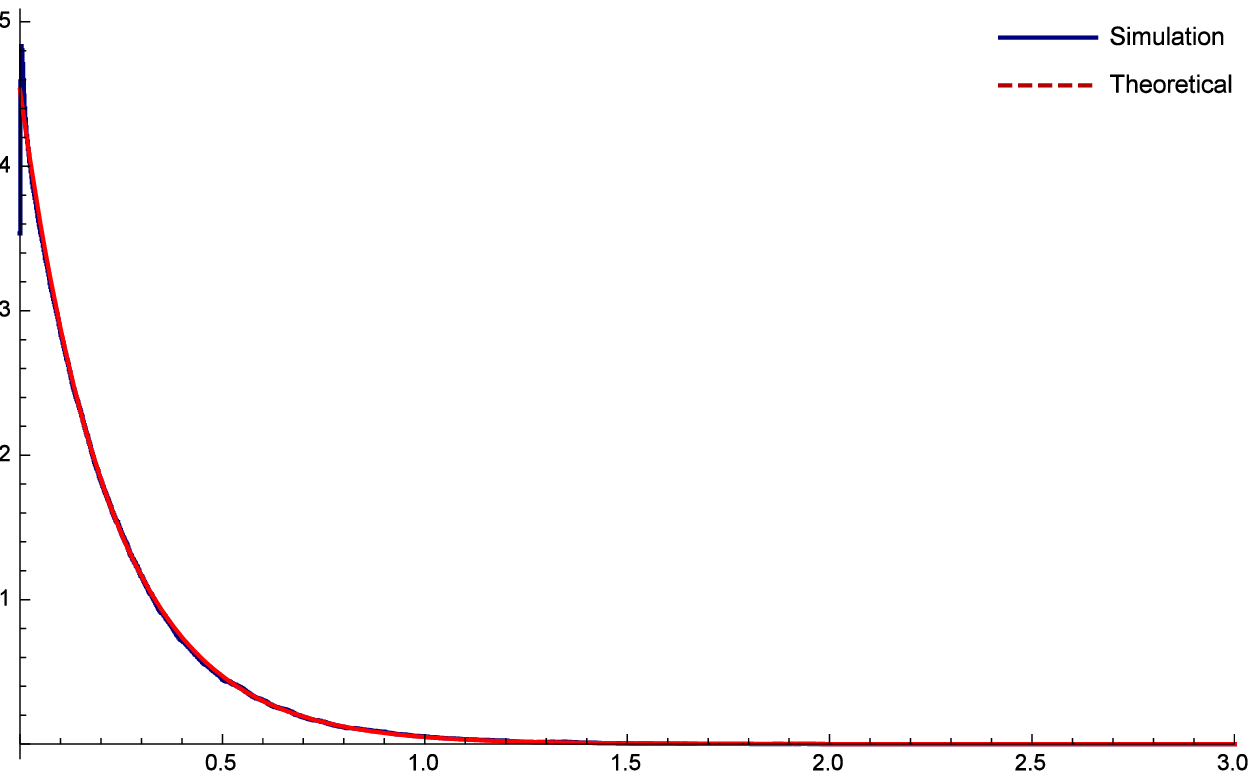}\\
(c) Queue 3 & (d) Queue 4 \\
\end{tabular}
\caption{\boldmath Simulated and theoretical scaled queue lengths for Example 4.}
\label{ql4b}\vspace*{-2ex}
\end{figure}

To obtain more insight in the influence of overloaded queues, we have chosen to do focus on one specific part of Figure \ref{fig:ql4}, namely the part between $\Lambda=\Lambda_3^\textit{crit}$ and $\Lambda=\Lambda_2^\textit{crit}$, i.e. the interval $[0.435, 0.577]$. In this interval, queues 3 and 4 are unstable, causing the system to behave like a two-queue polling system from the perspective of customers in queues 1 and 2. This corresponding system is similar to the (first two queues of the) original system, but has a longer switch-over time from $Q_2$ to $Q_1$, which is distributed as the convolution of the original switch-over times $S_2, S_3, S_4$, and $k_3$ services of type-3 customers, and $k_4$ services of type-4 customers. Figure \ref{ql4bcorresponding} depicts the mean queue lengths of queues 1 and 2 in the original model (taken from Figure \ref{fig:ql4}) and, additionally, in the corresponding two-queue polling model, while varying the total arrival rate $\Lambda$ between $0\leq \Lambda < 0.577$. It can clearly be seen in Figures \ref{fig:ql4}(a) and \ref{fig:ql4}(b) that the mean queue lengths of the two systems are different for $\Lambda < 0.435$, but they are the same in the aforementioned region $0.435 < \Lambda < 0.577$. This confirms that the two systems behave identically in this region. 

\begin{figure}[!htb]
\begin{tabular}{c}
\begin{tikzpicture}
\node[anchor=south west,inner sep=0] at (0,0) {\includegraphics[width=0.7\linewidth]{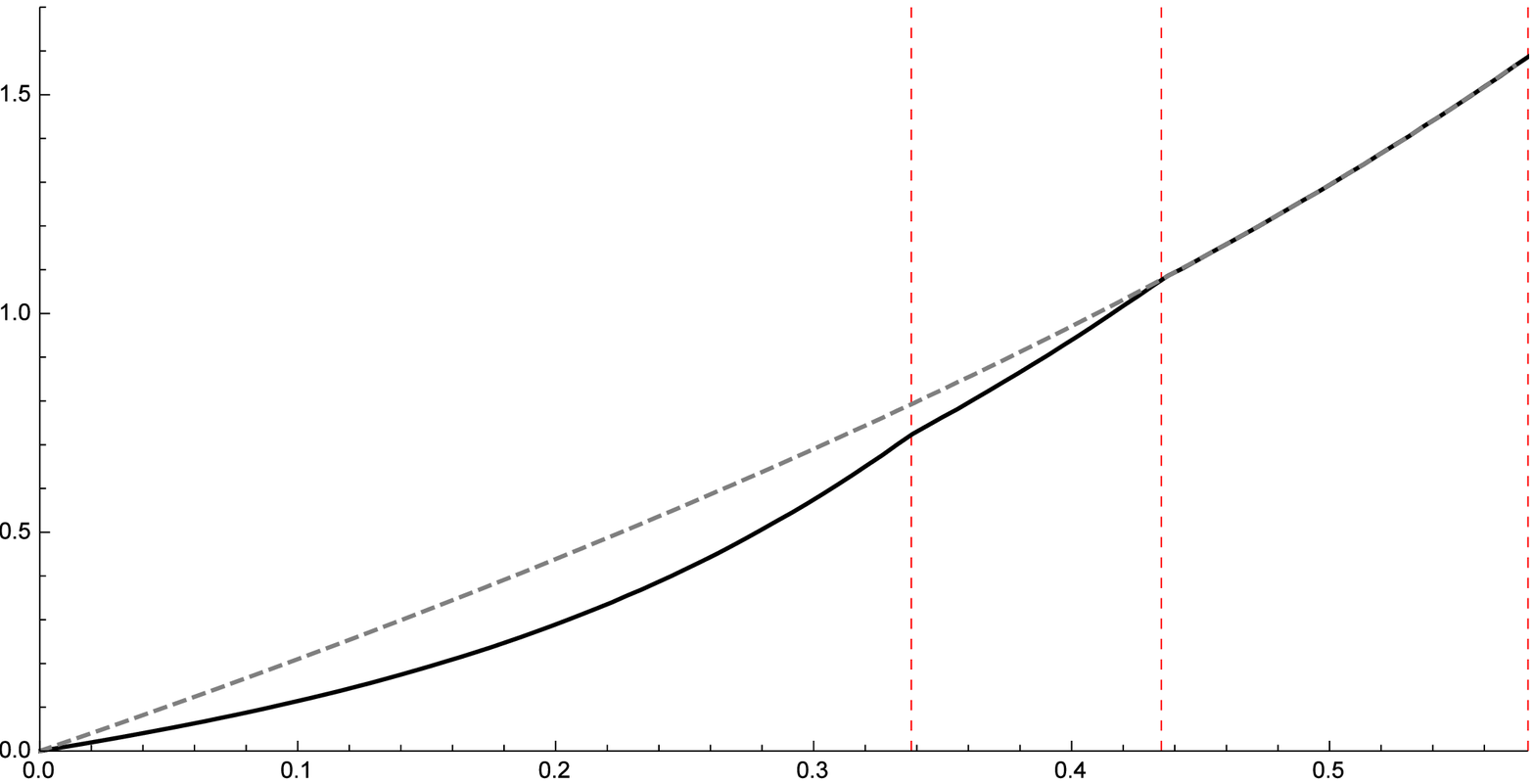} };
\node[anchor=south west,inner sep=0] at (-0.4,3.1) {\small $\E[N_1(\Lambda)]$};
\node[anchor=south west,inner sep=0] at (3.1,3.1) {\tiny $\Lambda_4^\textit{crit}$};
\node[anchor=south west,inner sep=0] at (4.2,3.1) {\tiny $\Lambda_3^\textit{crit}$};
\node[anchor=south west,inner sep=0] at (5.8,3.1) {\tiny $\Lambda_2^\textit{crit}$};
\node[anchor=south west,inner sep=0] at (2.9,-0.2) {\small $\Lambda$ };
\end{tikzpicture}
\\
(a) Queue 1\\
\begin{tikzpicture}
\node[anchor=south west,inner sep=0] at (0,0) {\includegraphics[width=0.7\linewidth]{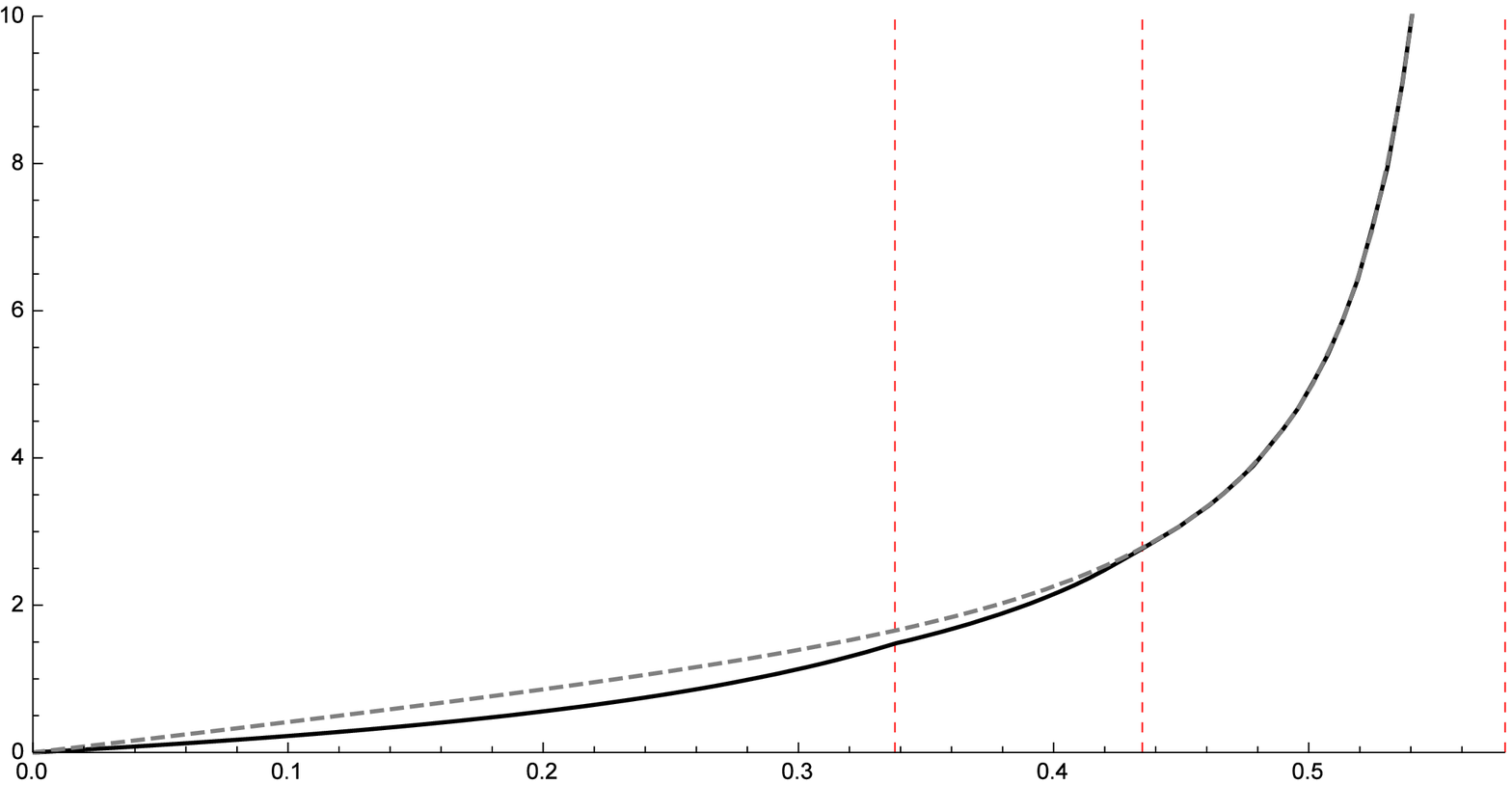} };
\node[anchor=south west,inner sep=0] at (-0.4,3.1) {\small $\E[N_2(\Lambda)]$};
\node[anchor=south west,inner sep=0] at (3.1,3.1) {\tiny $\Lambda_4^\textit{crit}$};
\node[anchor=south west,inner sep=0] at (4.2,3.1) {\tiny $\Lambda_3^\textit{crit}$};
\node[anchor=south west,inner sep=0] at (5.8,3.1) {\tiny $\Lambda_2^\textit{crit}$};
\node[anchor=south west,inner sep=0] at (2.9,-0.2) {\small $\Lambda$ };
\end{tikzpicture}
\\
(b) Queue 2\\
\end{tabular}
\caption{\boldmath Simulated mean queue lengths for queues 1 and 2 in the original model (solid lines) and in the ``corresponding'' polling model (dashed lines).}
\label{ql4bcorresponding}
\end{figure}

\vspace*{-3mm}
\section{Conclusion}

The present note has studied a novel application of the singular-perturbation technique leading to a rigorous proof of the HT limits in $k$-limited polling models with switch-over times.  The scaled queue-length of the critically loaded queue follows an exponentially distribution of which the parameter is known. The number of customers in the stable queue has the same distribution as the number of customers in a vacation system. Vacation systems with $k$-limited service can be analyzed very efficiently using modern numerical techniques. In particular, when the service times and switch-over times are constant, or in case they have a phase-type distribution or a Gamma distribution, the notorious root-finding problem becomes trivial  (see \cite{boonvanleeuwaarden2015,lee89} for more details). As such, the asymptotics form an excellent basis for approximating $k$-limited polling models - and, thus, WPAN and mobile (adhoc) networks
- with general load.

\section*{References}

\bibliographystyle{abbrv}

\end{document}